\documentclass[12pt]{article}
\usepackage{dblfloatfix}
\usepackage{amsmath, amssymb}
\usepackage{amsfonts}
\usepackage{times}
\usepackage{bm}
\usepackage[plain,noend]{algorithm2e}
\usepackage{color}
\usepackage{bm}
\RequirePackage[OT1]{fontenc}
\usepackage{verbatim,enumerate}
\usepackage{rotating, lscape}

\usepackage{graphicx,enumitem}
\newtheorem{problem}{Problem}
\newtheorem{prop}{Proposition}
\newtheorem{remark}{Remark}
\newtheorem{theorem}{Theorem}

\def\red{\textcolor{red}}
\def\R{\mathbb{R}}

\def\C{\mathbb{C}}

\newcommand{\N}{\mathbb{N}}
\newcommand{\Z}{\mathbb{Z}}
\newcommand{\F}{\mathfrak{F}}

\renewcommand{\Re}{\mathop{\rm Re}\nolimits}

\newcommand{\dmna}{_{\delta ,\mu ,\nu , \alpha}}

\begin{document}

\title{Buhmann covariance functions, their compact supports, and their smoothness}

\author{E. PORCU\footnote{Departamento de Matem\'atica, Universidad Federico Santa
Maria,}, V.P. ZASTAVNYI \footnote{Department of Mathematics, Donetsk National University}, and M. BEVILACQUA\footnote{Department of Statistics, University of Valparaiso.}}

\maketitle

\begin{abstract}
We consider the Buhmann class of compactly supported radial basis functions, whih includes a wealth of special cases that have been studied in both numerical analysis and spatial statistics literatures. In particular, the celebrated Wu, Wendland and Missing Wendland functions are notable special cases of this class. We propose a very simple difference operator
and show the conditions for which the application of it to Buhmann functions preserves positive definiteness on $m$-dimensional Euclidean spaces. We also show that the application of the difference operator increases smoothness at the origin,  whilst keeping positive definiteness in the same $m$-dimensional Euclidean space, as well as compact support. Thus, our operator is a competitor of the celebrated Mont{\'e}e operator, which allows to increase the smoothness at the origin, at the expense of losing positive definiteness in the space where the radial basis function is originally defined.
The proofs of our results highlight surprising connections with past literatures on celebrated class of functions. Amongst them, absolute and completely monotone functions.
\end{abstract}

%\begin{keywords}
%\end{keywords}

\def\hh{\boldsymbol{h}}
\section{Introduction}

Interpolation of data has a notable importance in both numerical analysis and geostatistical communities: in geostatistics, the underlying structure of the data is assumed to be a stochastic process, with the interpolation procedure known as kriging \cite{scheuerer}. The method is mathematically equivalent to kernel interpolation, a method used in numerical analysis for the same problem, but derived under completely different modelling assumptions.

Radial basis functions are well-known and successful tools for the
interpolation of data in many dimensions \cite{Buh}. Several radial basis functions of
compact support that give rise to nonsingular interpolation problems have
been proposed, and we cite \cite{euclid1, euclid2, Buh, Wen, Schaback, wu, Zast2000} and the impressive reviews in \cite{fassauer} and \cite{schaback-wendland}, amongst others. In particular, the motivation behind Buhmann's \cite{Buh} tour de force is to propose radial functions of compact support
that also give positive definite matrices and have genuinely banded interpolation
matrices (similarly to the multiquadric and Gaussian kernels). Of such nature are those ones we will discuss in the present paper. Early examples of radial functions with compact support that have a simple
piecewise polynomial structure are due to \cite{wu}. Then Schaback and Wu \cite{s-w},
 Wendland \cite{Wen} and finally Schaback \cite{Schaback} established several of their special properties,
such as certain optimality facts about their degree and smoothness. The functions proposed by Buhmann are closely related to the so-called multiply monotone radial
basis functions as discussed in \cite{micchelli}.

Radial basis functions are known under the name of covariance (correlation) functions in the geostatistical community: the use of compactly supported covariance functions has been advocated in a number of papers, and we refer the reader to \cite{gnei02}, with the references therein, and to \cite{DPB14} for a recent effort under the framework on multivariate Gaussian fields. Covariance functions with compact support
represent the building block for the construction of methods allowing to overcome the big data problem (\cite{furrer}). The recent work of \cite{BFFP} brought even more attention on the role of some classes of compactly supported covariances for asymptotically optimal prediction on a bounded set of $\R^d$. From the cited works it has become apparent that the smoothness at the origin (intended as even extension) of a compactly supported and isotropic covariance function plays a crucial role for both estimation and prediction. Wendland functions \cite{Wen} have been especially popular, being compactly supported
over balls of $\R^d$ with arbitrary radii, and additionally allowing for
a continuous parameterization of differentiability at the origin, in a
similar way to the Mat{\'e}rn family (\cite{stein}).

\def\red{\textcolor{red}}

The tour de force in \cite{Wen}, \cite{gnei02} and \cite{zucastell1} puts emphasis on linear operators, called {\em Mont{\'e}e}, that allow to increase the smoothness of a given radial function, being positive definite on $\R^d$. This is done at the expense of losing positive definiteness, which is only achieved on $\R^{d-2}$, for $d \ge 3$.

Figure 1 depicts the following situation: Dashed lines report Wendland functions \cite{Wen} with unit compact support, for $k=0,1,2$ (from left to right) being the parameter that allows to govern differentiability at the origin. Dashed-dotted lines report the respective Wendland functions for a compact support equal to $0.75$. The functions depicted with continuous lines report their weighted differences (see Equation (\ref{del_piero})). We can clearly appreciate that the level of differentiability changes for these last.

This gives a substantial motivation for the present work: we
consider the Buhmann functions as defined in
Zastavnyi~\cite{zast2006} on the basis of \cite{Buh}. This class
includes as special case many other classes of compactly supported
covariance functions, such as Askey \cite{askey} Wendland
\cite{Wen} and Missing Wendland \cite{Schaback} functions,  as well as the  Zastavnyi \cite{c11, c10,
zast2006} and Trigub \cite{Trigub, TrBel} classes. Finally, also Wu
functions \cite{wu} and the celebrated spherical model \cite{wack}
are included as special cases. We show that some parameterized
differences of Buhmann functions preserve positive definiteness in
$m$-dimensional Euclidean spaces, and we then determine the exact
level of smoothness induced by such operation.

\def\bx{\boldsymbol{x}}

%Let $\{Z(\bx), \bx \in D \}$ be a stationary and isotropic Gaussian
%field observed over a bounded set $D$ of $\R^m$. The assumption of
%Gaussianity implies that we only need focus on the first and second
%order moments in order to specify the probabilistic properties of
%the field.
We focus throughout on the class $\Phi_m$
of continuous functions $\varphi:[0,\infty) \to \R$  such that
$\varphi(\|\bx\|)$ is positive definite on $\R^m$. Thus,
$\varphi(\|\cdot\|)$ with $\varphi(0)=1$ is the correlation function
of some Gaussian field.

Apparently, the
functions $C(\cdot):=\varphi(\| \cdot \|)$ are radially symmetric and Schoenberg's theorem
(1938, see \cite{daley-porcu} for a more recent discussion) uniquely
identifies them as scale mixtures of the type
\begin{equation}
 \varphi(t)= \int_{[0,\infty)} \Omega_m(r t) F({\rm d} r), t \ge 0,
 \end{equation}
with $F$ a uniquely determined probability measure, and $\Omega_m(\cdot)$ being the characteristic function of a random vector that is uniformly distributed on the spherical shell of $\R^m$. Daryl and Porcu \cite{daley-porcu} put emphasis on the measure $F$, termed Schoenberg measure there. It is well known \cite{gnei02} that any random vector $\boldsymbol{X}$ of $\R^m$ with characteristic function $\varphi$ can be written as $\boldsymbol{X}= \boldsymbol{\eta} R$, with $\boldsymbol{\eta}$ having $\Omega_m$ as characteristic function, $R$ a positive random variable distributed according to $F$, and $\boldsymbol{\eta}$ and $R$ are stochastically independent. The identity above is intended as equality in distribution. \\
The class $\Phi_m$ is nested, with the following inclusion relation
$$ \Phi_1 \supset \Phi_2 \supset \ldots \supset \Phi_{\infty}:= \bigcap_{m \ge 1} \Phi_m, \qquad m \in \mathbb{N}, $$
being strict (see, for example, \cite{daley-porcu,gnei02}).
%This will be one the key arguments to prove the solutions to the problems described subsequently.
A function $f:(0,\infty) \to \R$ is called completely monotone if it
is infinitely often differentiable and $(-1)^n f^{(n)}(x) \ge 0$,
for all $n \in \mathbb{Z_+}$ and for all $x>0$. The set of
completely monotone functions on $(0,\infty)$ is denoted ${\cal
CM}$. By Schoenberg's theorem, $\varphi(t) \in \Phi_{\infty}$ if and
only if $\varphi(\sqrt{t}) \in {\cal CM}$ with
$\varphi(0)=\varphi(+0)<\infty$.

The plan of the paper is the following. Section 2 introduces the Buhmann class and offers a formal statement of the problem. Section 3 exposes the structure of the solution. We also characterize the smoothness of differences of Buhmann functions in a neighborhood of the origin. Section 4 illustrates connections with previous literature.

\section{Buhmann functions. Statement of the problem}

\subsection{Buhmann's class and relations with previous literature}

We denote $C(\R^m)$ the set of continuous functions on $\R^m$, for $m=1,2,\ldots.$
 Let $\delta , \mu , \nu \in \C _+:=\left\{z\in \C:\Re z>0 \right\}$ and $\alpha\in\C$.
 {Zastavnyi~\cite{zast2006} (2006) proposed the following even functions given on $\R$:}
   %%%
\begin{equation} \label{bu}
\varphi_{\delta,\mu,\nu,\alpha}(x):=\left\{
\begin{array}{ll}
\int^{1}_{|x|}  (s^2-x^2)^{\nu-1}(1-s^{\delta})^{\mu-1}
s^{\alpha-2\nu+1}
 \,{\rm d} s ,& \qquad  |x|<1\\
 & \\
0, & \qquad |x|\ge 1.
\end{array}\right.
\end{equation}

If $\delta,\mu,\nu\in\C_+$, then {arguments in Proposition 1 and Theorem 1 in
\cite{zast2006} show, respectively, that} $\varphi_{\delta,\mu,\nu,\alpha}\in
C(-1,1)$ {if and only if} $\alpha\in\C_+$ %(see \cite[Proposition~1]{zast2006})
and that $\varphi_{\delta,\mu,\nu,\alpha}\in C(\R)$ {if and only if}
$\alpha,\mu+\nu-1\in\C_+$  If $\delta,\mu,\nu,\alpha\in\C_+$, then
$\varphi_{\delta,\mu,\nu,\alpha}(0)=B(\alpha/\delta,\mu)/\delta$, {with $B$
denoting the Beta function.}

The functions $\varphi_{\delta,\mu,\nu,\alpha}$ coincide (modulo some positive factors) with the
functions $$\phi_{\delta,\varrho, \lambda, \alpha}(x)\equiv
 2\varphi_{2\delta,\varrho+1,\lambda+1,
2\alpha+2}(x), \qquad x \in \R, $$ introduced by Martin Buhmann \cite{Buh}. We thus term them Buhmann functions throughout.
The class $\varphi_{\delta,\mu,\nu,\alpha}$ includes a wealth of interesting special cases. For instance,
 $\mu\delta
\varphi_{\delta,\mu,1,\delta}(x)=(1-|x|^{\delta})^{\mu}_{+}$,  {which implies
that $\mu \varphi_{1,\mu,1,1}$ coincides with the Askey functions
\cite{askey}.} Also, we have that
%%%
%In special cases, we get the following known functions:
\begin{equation}
\varphi_{1,\mu,\nu,2\nu -1}(x) \equiv h_{\mu,\nu}(x)\equiv
\frac{2^{\nu-1}\Gamma(\nu)}{\mu} \psi_{\mu,\nu-1}(x), \qquad x \in \R,
\end{equation}
with the functions $h_{\mu,\nu}$ being introduced by Zastavnyi
(2002) \cite{c10,c11} and defined as follows: $h_{\mu,\nu}(x):=0$ for $|x|\ge 1$
and
\begin{equation} \label{f38}
\begin{split}
h_{\mu,\nu}(x):=&
 \int_{|x|}^{1} (2u-|x|)g_{\mu,\nu}(u)g_{\mu,\nu}(u-|x|)\;du,\;|x|<1,
 \\ \text{ where } g_{\mu,\nu}(u):=&
 u^{\mu-1} (1-u^2)^{\nu-1},\, u\in (0,1),\;\mu,\nu\in\C_+ .
\end{split}
\end{equation}
%%%%%%
Functions of the form \eqref{f38} arise in the study of exponential
 type entire functions without zeros in the lower half-plane~\cite[Proposition~5.1]{Zast2006_UMB}.
 The functions $\psi_{\mu,\nu-1}$, with
$\mu>0$, $\nu\in\N$, have been introduced in 1995 by Wendland \cite{Wen}, {and
they have been termed Wendland functions in both numerical analysis and
geostatistical literatures}: for $\mu>0, {k\in\Z_{+}}$, we have
\begin{equation*}
\psi_{\mu,0}(x):=\psi_{\mu}(x):=(1-|x|)^{\mu}_{+},\; \;
\psi_{\mu,k}:=I^k \psi_{\mu}\;({k\in\N}),
\end{equation*}
where { $I(f)(x):=\int_{|x|}^{+\infty} sf(s)\,{\rm d} s$ is the
Matheron's \cite{mat} Mont{\'e}e operator (provided the integral is
well defined), and where $I^k$ is the $k$-fold application of the
operator $I$.}
%Walks through dimension became especially popular after Gneiting's seminal paper \cite{gnei02}.
Arguments in \cite{Wen} and subsequently \cite{gnei02} show that $I
\varphi$ belongs to the class $\Phi_{m-2}$ whenever $\varphi \in
\Phi_m$, for $m \ge 3$. For $k<2m$, the $k$-fold application of the
Mont{\'e}e operators shows that ${I}^k \varphi \in \Phi_{m-2k}$,
$k\in\N$.

Gneiting  \cite[Equation (17)]{gnei02} has proposed a
generalization of Wendland functions on the basis of the fractional
Mont{\'e}e operator, which coincides with the {normalized Buhmann}
functions
 $\varphi_{1,\mu+1,\nu,2\nu}(x)/\varphi_{1,\mu+1,\nu,2\nu}(0)$, $\mu, \nu>0$,
  as well as with the functions $h_{\mu,\nu+1}(x)/h_{\mu,\nu+1}(0)\equiv \psi_{\mu,\nu}(x)/\psi_{\mu,\nu}(0)$ (see Equation~\eqref{eq_all}).
% and the functions $\psi_{\mu,\nu-1}$, with
%$\mu,\nu\in\N$, being proposed by Wendland (1995) \cite{Wen} and Gneiting \cite{gnei02}.
 Arguments in \cite{Buh} show
that Wu functions \cite{wu} and consequently the spherical model are {special
cases} of the Buhmann class.

 For $r\in\Z_+$ and $k\in\N$, we
have $$h_{r+k,r+1}(x)\equiv B(r+k,2r+1) A_{r,2k-1}(x),$$ with the
splines $A_{r,2k-1}$ introduced by Trigub  (1987), and we refer to
\cite{Trigub}, \cite[\S~6.2.13, 6.2.16, 6.3.12]{TrBel} for their
analytical expression which is not reported here. {Equation above in
turn highlights the explicit connection between Trigub splines and
Wendland functions:
$A_{r,2k-1}(x)\equiv\psi_{r+k,r}(x)/\psi_{r+k,r}(0)$, for $r \in
\Z_+$ and $k \in \N$.}

For a proof of the identities
above, the reader is referred to Zastavnyi and Trigub~\cite[Remarks
10 and 11]{c10}, to \cite[Theorems
12 and 13]{c11}, \cite{zast2006} and~\cite[\S~4.7]{zastavnyi:2008}. \\
%{\bf I think it is more accurate to put a reference close to each identity}.
Arguments in Proposition 4 of \cite{zast2006} show that, for $\delta,
\mu,\nu\in\C_+$ and  $x\in\R,$ \\

\begin{equation}
\begin{split}
&
 \varphi_{2,\frac{\mu}{2},\frac{\mu}{2}+\nu,2\nu-1}(x)\equiv \frac{2^{\mu-1}\Gamma\left(\frac{\mu}{2}\right)\Gamma\left(\frac{\mu}{2}+\nu\right)}{\Gamma(\mu)}
 \varphi_{1,\mu,\nu,2\nu-1}(x)\,,\;
%\mu,\nu\in\C_+,\; x\in\R,
\\&
 2\nu\varphi_{\delta,\mu+1,\nu,2\nu}(x)\equiv\delta\mu
 \varphi_{\delta,\mu,\nu+1,2\nu+\delta}(x)\,,\;
 %\delta, \mu,\nu\in\C_+,\; x\in\R,
\end{split}
\end{equation}
\\
and, for {$\mu,\nu\in\C_+$ and $x\in\R$,} we also have the obvious
identities:
%Obviously
\begin{equation}\label{eq_all}
 \varphi_{1,\mu+1,\nu,2\nu}(x)
{\equiv}
 \frac{\mu}{2\nu}\varphi_{1,\mu,\nu+1,2\nu+1}(x)\equiv
  \frac{\mu}{2\nu}h_{\mu,\nu+1}(x)\equiv
   2^{\nu-1}\Gamma(\nu)\psi_{\mu,\nu}(x).
\end{equation}

\subsection{{Buhmann class and its Fourier and Laplace Transforms}}

{After the illustration of the relation between Buhmann and other celebrated
classes of radial basis function, we need some preliminary material in order to
provide a better description of the results coming subsequently.}
%Let us now recall some preliminary material
For {a function} $h$ {defined on} $(0,\infty)$ and  $m\in\C$, we
define the Hankel transform $\F_m$ as follows:
\begin{equation} \label{hankel}
   \F_m(h)(t):= t^{1-\frac{m}{2}}
             \int_{0}^{\infty} h(u) u^{\frac{m}{2}} J_{\frac{m}{2}-1} (tu)\,{\rm d}u=
             \int_{0}^{\infty} h(u) u^{m-1} j_{\frac{m}{2}-1}
             (tu)\,{\rm d}u, \qquad t >0,
\end{equation}
where   $J_{\lambda}$ is the Bessel function of the first kind (see
\cite[Sec.~3.1]{Watson})  and
\begin{equation}\label{bes}
j_{\lambda}(x):=\frac{J_{\lambda}(x)}{x^{\lambda}}
 =\frac{1}{2^{\lambda}} \sum_{k=0}^{\infty}
 \frac{1}{\Gamma(k+\lambda+1)}\cdot
 \frac{\left(-\frac{x^2}{4}\right)^{k}}{k!}\;,\;x\in\C\;,\;\lambda\in\C\,.
\end{equation}
%%%%%%%%%%%%%%%%%%
 \begin{remark}\label{rem_1}
 {\rm For $m\in\N$ the transform $\F_m$ is connected with the Fourier
transforms $F_m$ of radial functions through the identity
 $$
 F_m(h(\|\cdot\|))(\bx)=(2\pi)^{\frac{m}{2}}
 \F_m(h)(\|\bx \|) \,, \, \qquad \bx\in\R^m\,.
 $$
These facts and Bochner-Khintchine theorem (see, for example
\cite{Sasvari_2013, Schilling, TrBel}) imply that, if $h$ is a continuous
functions on $[0,\infty)$ and $\int_{0}^{\infty}t^{m-1}|h(t)|\,dt<\infty$, then
$h \in \Phi_m$ if and only if $\F_m(h)$ is nonnegative on the positive real
line.}
 \end{remark}

For $\delta , \mu ,  \alpha+1\in \C_+ $ and $\nu\in\C$, we define
the function $I_{\delta,\mu,\nu,\alpha}: \R_+ \to \C$ through
%  and the function $I_{\delta,\mu,\nu,\alpha}$ for   $\delta , \mu ,  \alpha+1\in \C_+ $ and $\nu\in\C$ we define the formula
%%
 \begin{eqnarray} \label{I}
  I_{\delta,\mu,\nu,\alpha}(t)&:=&  t^{-\alpha-1-\delta
(\mu-1)} \int_{0}^{t}
(t^{\delta}-u^{\delta})^{\mu-1}u^{\alpha-\nu+\frac{1}{2}}
J_{\nu-\frac{1}{2}} (u) \;{\rm d}u \nonumber \\ &=& \int_{0}^{1}
(1-x^{\delta})^{\mu-1}x^{\alpha} j_{\nu-\frac{1}{2}} (tx) {\rm d} x
  \,,\,\qquad \qquad t>0 \,.
 \end{eqnarray}

The following result reports succinctly a collection of useful
results from~\cite{zast2006} (Theorems~2,~3 and
Proposition~4~(Assertions~1,3)).

\begin{theorem}[Zastavnyi \cite{zast2006}] \label{leo_messi} Let the functions $I$ and $\F_m$ as being defined through Equations (\ref{I}) and (\ref{hankel}), respectively. Denote with $I'$ the first derivative of $I$. Then, the following assertions are true:
\begin{enumerate}
%{\rm \cite[Theorem~3]{zast2006}}.
\item  Let $\delta , \mu , \nu ,  m , \alpha+m\in \C_+ $. Then
$\F_m(\varphi_{\delta,\mu,\nu,\alpha})(t)=2^{\nu-1}\Gamma(\nu)
I_{\delta,\mu,\frac{m-1}{2}+\nu,m-1+\alpha}(t).$ Moreover, if
 $n , m-n+2\nu \in \C_+ $, then
 $$\F_m(\varphi_{\delta,\mu,\nu,\alpha})(t) =
 \frac{2^{\frac{n-m}{2}}\Gamma(\nu)}{\Gamma(\frac{m-n}{2}+\nu)}
 \F_n(\varphi_{\delta,\mu,\frac{m-n}{2}+\nu,m-n+\alpha})(t), \qquad t \ge 0.$$
 % {\rm \cite[Proposition~4
%(1)]{zast2006}}
 \item Let $\delta,\mu,\alpha+1\in\C_+$ and
$\nu\in\C\,$. Then
\begin{equation}\label{42}
\begin{split}
  I'\dmna(t)&=-tI_{\delta, \mu, \nu+1,\alpha+2}(t), \qquad t>0\,.
\end{split}
\end{equation}
%  {\rm \cite[Proposition~4 (3)]{zast2006}}.
\item If
$\mu,\nu\in\C_+$, $t>0$, then
 \begin{equation}\label{44}
 I_{1,\mu,\nu,2\nu-1}(t)=\frac{2^{\frac{1}{2}-\nu}\Gamma(\mu)\Gamma(2\nu)}{\Gamma\left(\nu+\frac{1}{2}\right)\Gamma(\mu+2\nu)}
 \,{_1}F_2\left(\nu ; \frac{\mu+2\nu}{2} , \frac{\mu+2\nu+1}{2} ;
 -\frac{t^2}{4}\right)\,.
\end{equation}
%{\rm \cite[Theorem~2 (1)]{zast2006}}.
\item   Let $\delta ,
\mu , \nu \in\C_+$ and $\alpha \in\C$. Then
%\textcolor{red}{
 \begin{equation}\label{rek}
\varphi_{\delta,\mu,\nu+1,\alpha+2}(t)=2\nu\int_{|t|}^{\infty}
u\,\varphi_{\delta,\mu,\nu,\alpha}(u)\,{\rm d}u\;,\;t\ne 0\,.
 \end{equation}
% {\bf Viktor: $t$ or $x$? you used $x$ in Equation (\ref{bu}).}
\end{enumerate}
\end{theorem}

A relevant remark is that Equation (\ref{44}) describes the spectral density of the Wendland functions.
Another remarkable consequence of Theorem \ref{leo_messi} is that, for $\mu,\nu,m,2\nu-1+m\in\C_+$,
\begin{equation}\label{F_m}
\F_m(h_{\mu,\nu})(t)=\F_m(\varphi_{1,\mu,\nu,2\nu-1})(t)=
 %% \frac{2^{\frac{1-m}{2}}\Gamma(\nu)}{\Gamma(\frac{m-1}{2}+\nu)}   \F_1(\varphi_{1,\mu,\frac{m-1}{2}+\nu,m-1+2\nu-1})(t)=
 \frac{2^{\frac{1-m}{2}}\Gamma(\nu)}{\Gamma(\frac{m-1}{2}+\nu)}
\F_1(h_{\mu,\frac{m-1}{2}+\nu})(t)\,, \qquad t \ge 0,
\end{equation}
which in turn shows, in concert with \cite[Lemma~12]{c10}, that in some cases the Hankel transforms above can be written in closed form. Specifically, we have
\begin{equation}\label{Fh}
\begin{split}
&\F_m(h_{\mu,\nu})(t)=\F_m(\varphi_{1,\mu,\nu,2\nu-1})(t) =
2^{\nu-1}\Gamma(\nu) I_{1,\mu,\frac{m-1}{2}+\nu,m-1+2\nu-1}(t)=
\\&
  D(m,\mu,\nu)\cdot
 {_1}F_2  \left(\frac{m-1}{2}+\nu;\frac{m-1}{2}+\nu+\frac{\mu}{2} ,\frac{m-1}{2}+\nu+\frac{\mu+1}{2} ;-\frac{t^2}{4} \right)
\,,
\\&
\text{with }
 D(m,\mu,\nu):=   \frac{2^{-\frac{m}{2}}\Gamma(\nu)\Gamma(\mu)\Gamma(m-1+2\nu)}{\Gamma\left(\frac{m}{2}+\nu\right)\Gamma(\mu+m-1+2\nu)},\;  \mu,\nu,m,2\nu-1+m\in\C_+\,.
 \end{split}
\end{equation}

Let us use the abuse of notation $\widehat{h}_{\mu,\nu}$ for the one dimensional Fourier transform of the function $h_{\mu,\nu}$. We also denote with $L$ the Laplace transform operator.
For $\mu,\nu\in\C_+$, arguments in Zastavnyi and Trigub \cite[Equation (44)]{c10} show that
\begin{equation}\label{f44}
L\left( t^{2\nu+\mu -1} \widehat{h}_{\mu,\nu}(t) \right) (x):=
\int_{0}^{\infty}e^{-tx}\, t^{2\nu+\mu -1}
\widehat{h}_{\mu,\nu}(t)\,{\rm d}t=
 \frac{\Gamma ^2 (\nu )  \Gamma (\mu) 2^{2\nu
-1}}{x^\mu (1+x^2)^\nu },\,  x>0.
\end{equation}
Thus, for $ \mu,\nu,m,2\nu-1+m\in\C_+$ and $x>0$, we have
\begin{equation}\label{pdh}
\begin{split}
& L\left( t^{m-1+2\nu+\mu -1} \F_m(h_{\mu,\nu})(t) \right) (x)
 \\&=
\frac{2^{\frac{1-m}{2}}\Gamma(\nu)}{\Gamma(\frac{m-1}{2}+\nu)}
  L\left( t^{m-1+2\nu+\mu -1} \F_1(h_{\mu,\frac{m-1}{2}+\nu})(t) \right) (x)
  \\&=
  \frac{2^{\frac{1-m}{2}}\Gamma(\nu)}{\Gamma(\frac{m-1}{2}+\nu)}\cdot
  \frac{\Gamma ^2 ({\frac{m-1}{2}+\nu} )  \Gamma (\mu) 2^{m-1+2\nu -1}}{(2\pi)^{\frac{1}{2}}}\cdot
  \frac{1}{x^\mu (1+x^2)^{\frac{m-1}{2}+\nu} }
 \\&=
 \frac{C(m,\mu,\nu)}{x^\mu (1+x^2)^{\frac{m-1}{2}+\nu} } \, ,
  \\&
 \text{with } C(m,\mu,\nu):=
  \frac{2^{-\frac{m}{2}}\Gamma(\nu)\Gamma(\mu)\Gamma(m-1+2\nu)}{\Gamma\left(\frac{m}{2}+\nu\right)}
  =D(m,\mu,\nu)\Gamma(\mu+m-1+2\nu)\,.
\end{split}
\end{equation}

\begin{remark}\label{rem_2}
{\rm Equation~\eqref{pdh}  is the crux of the proof of the main part of
Theorem~11 in~\cite{c11}:}\\
 {\it
$(i)$~If $\nu>\frac12$ and $\mu\ge\max\{\nu,1\}$, then
$h_{\mu,\nu}\in \Phi_1$. If, additionally, $(\mu;\nu)\ne (1;1)$,
then there exist constants $c_i>0$, {$i=1,2$}, depending on $\mu$
and $\nu$ only, such that $$c_1\le (1+t^2)^{\nu}
\cdot\widehat{h}_{\mu,\nu}(t) \le c_2, \qquad t\in\R.$$
\\
$(ii)$ If $\nu \ge 1$, then $h_{\mu,\nu}\in \Phi_1 \iff \mu\ge\nu$.
\\
$(iii)$ If $m\ge 2$, then $h_{\mu,\nu}\in\Phi_m \iff \nu>\frac12$
and $\mu\ge \frac{m-1}{2}+\nu$. In this case, there exist two
constants $c_i>0$, {$i=1,2$}, depending on $\mu$, $\nu$ and $m$, and
such that
  $$c_1\le
(1+t^2)^{\frac{m-1}{2}+\nu} \cdot\F_m(h_{\mu,\nu})(t) \le c_2,\;
\qquad t\ge 0.$$ }
 {\rm This theorem is related to the positiveness of the
function $I_{1,\mu,\nu,2\nu-1}(t)$ for all $t>0$. Theorems on positiveness of
the functions  $I_{\delta,\mu,\nu,\alpha}(t)$ are obtained in~\cite[Theorems
4,5,6]{zast2006} (the well-known cases given before Theorem 4
from~\cite{zast2006}).}
\end{remark}

\subsection{Statement of the Problem}

We are now able to state our main
\begin{problem} \label{leo_messi2}
\em
 Let $\mu>0$, $\nu>\frac{1}{2}$,
$\mu+\nu>1$.
 Then  $h_{\mu,\nu}\in C(\R)$ (see~\cite{c10,c11}). Let
 $\varepsilon>0$ and
\begin{equation} \label{del_piero}
 f_{\mu,\nu,\varepsilon,\beta_1,\beta_2}(x):=
\beta_2^\varepsilon\, h_{\mu,\nu}\left(\frac{x}{\beta_2}\right)-
\beta_1^\varepsilon\, h_{\mu,\nu}\left(\frac{x}{\beta_1}\right),
\;x\in\R.
 \end{equation}
Let $m\in\N$. Show the conditions on $(\varepsilon, \mu, \nu)$ such
that, for any $\beta_2>\beta_1>0$, we have
\begin{equation}\label{strong}
f_{\mu,\nu,\varepsilon,\beta_1,\beta_2}\in\Phi_m\,.
\end{equation}
 \end{problem}

Next section details the structure of the solution and offers the
exact smoothness of the new covariance resulting from the
differences of Buhmann functions.

%\newpage
\section{Structure of the solution}

Let us start with a general assertion regarding the structure of Problem \ref{leo_messi2}.

\begin{prop} \label{leo_messi3}
The following conditions are equivalent:
\begin{enumerate}
\item Condition \eqref{strong} is satisfied.
\item For any $\beta_2>\beta_1>0$, the function
  $t \mapsto \beta_2^{\varepsilon+m}\F_m(h_{\mu,\nu})(\beta_2t)-\beta_1^{\varepsilon+m}\F_m(h_{\mu,\nu})(\beta_1t)$ is nonnegative in interval $(0,\infty)$.
\item The function
$t^{\varepsilon+m}\F_m(h_{\mu,\nu})(t)=2^{\nu-1}\Gamma(\nu)t^{\varepsilon+m}I_{1,\mu,\frac{m-1}{2}+\nu,m-1+2\nu-1}(t)$
increases in the interval $(0,\infty)$.
\item The following inequality is true: \begin{equation*}
  \begin{split}
  &(\varepsilon+m)I_{1,\mu,\frac{m-1}{2}+\nu,m-1+2\nu-1}(t)+tI^{\,\prime}_{1,\mu,\frac{m-1}{2}+\nu,m-1+2\nu-1}(t)=
  \\&
  (\varepsilon+m)I_{1,\mu,\frac{m-1}{2}+\nu,m-1+2\nu-1}(t)-t^2I_{1,\mu,\frac{m-1}{2}+\nu+1,m-1+2\nu+1}(t)\ge
  0,\,\forall t>0
  \end{split}
  \end{equation*}
\item Let $n:=\frac{m-1}{2}+\nu$. Then,
$$(\varepsilon+m)\F_1(h_{\mu,n})(t)-\frac{t^2}{2{n}}\F_1(h_{\mu,n+1})(t)\ge
0,\, \qquad \qquad \forall t>0\,.$$
\item We have
 \begin{equation*}
  \begin{split}
  &
L\left( t^{2n+\mu -1}\Big( (\varepsilon+m)\widehat{h}_{\mu,n}(t)-
 \frac{t^2}{2{n}}\widehat{h}_{\mu,n+1}(t)\Big)\right) (x)=
 \\&
 \Gamma ^2 (n )  \Gamma
(\mu) 2^{2n -1}\left( \frac{\varepsilon+m}{x^\mu (1+x^2)^n }-
\frac{2n}{x^\mu (1+x^2)^{n+1} }\right)\in {\cal CM}\,.
\end{split}
  \end{equation*}
\item
 \begin{equation}\label{answer}
 \frac{\varepsilon-2\nu+1+(\varepsilon+m)x^2}{x^\mu (1+x^2)^{\frac{m-1}{2}+\nu+1}}\in {\cal CM}
 \end{equation}
 \end{enumerate}
 \end{prop}

%%%%%%%%%%%%%%
 \begin{remark}\label{re5}
 {\rm It follows from Hausdorff-Bernstein-Widder theorem  (see, for example, \cite{Feller_39, Sasvari_2013, Schilling, Widder}) that if
$g\in C{[0,+\infty)}$ and its Laplace transform
$$
Lg(x):=\int _0^{+\infty}e^{-xs}g(s)\ ds
$$
converges for all $x>0$, then $g(s)\ge 0$ for $s\ge 0$ if and only if $Lg\in
{\mathcal CM}$. }
\end{remark}
%%%%%%%%%%%%%%%%%
 The proof Proposition~\ref{leo_messi3} is an easy consequence of Remark~\ref{rem_1} in concert with the
Hausdorff-Bernstein-Widder theorem (see Remark~\ref{re5}),
Theorem~\ref{leo_messi} (statements 1 and 2), and equalities~\eqref{F_m}
and~\eqref{f44}.

 Note that, if $\mu,\nu >0$, then $x^{-\mu} (1+x^2)^{-\nu}\in
 {\cal CM}$ if and only if $\widehat{h}_{\mu,\nu}(t)\ge 0$ for all
 $t>0$,  if and only if $I_{1,\mu,\nu,2\nu-1}(t)\ge 0$ for all
 $t>0$ (see Equations~\eqref{f44} and~\eqref{Fh}).

Proposition \ref{leo_messi3} will now be combined with the following facts:

\begin{enumerate}
\item If $\nu\ge 1$, then $x^{-\mu} (1+x^2)^{-\nu}\in
 {\cal CM}$ if and only if $\mu\ge\nu$. The sufficiency of this result can be found in \cite{fields}, and the necessity has been proved in \cite{moak}, \cite[Lemma~8]{Zast2000}).
 \item  If $0<\nu<1$, $\mu\ge 1$, then $x^{-\mu} (1+x^2)^{-\nu}\in
 {\cal CM}$ \cite{moak} and \cite[Example 5.4]{Zast2006_UMB}, \cite[\S~4.7, Example 4.7.7]{zastavnyi:2008}.
\item If $\nu>0$, $\mu\ge 2\nu$, then $x^{-\mu} (1+x^2)^{-\nu}\in
 {\cal CM}$ \cite{askey-pollard}.
\item If $n=1,2,3$, then $(a+x^2)/(x^n (1+x^2)^n)\in
 {\cal CM}$ if and only if $ a\ge 1/(2^{n-1}+1)$ \cite[\S~2]{c10}.
 \end{enumerate}

We can now combine the first three sufficient conditions above to obtain the following assertion:
if $\nu>0$ and
 $\mu\ge\min\{2\nu;\max\{1,\nu\}\}$ , then $x^{-\mu} (1+x^2)^{-\nu}\in
 {\cal CM}$.

The combination of these facts with Proposition \ref{leo_messi3} has just offered the proof of the following
\begin{theorem} \label{cr7} The following assertions are true:
\begin{enumerate}
\item
 If \eqref{strong} is true, then $\varepsilon\ge
 2\nu-1$.
 \item
If $m\in\N$, $\nu>\frac{1}{2}$, $\varepsilon\ge 2\nu-1$ and
 $\mu\ge {(m-1)/2}+\nu+3$, then condition \eqref{strong} is
 true. If, in addition, $\varepsilon= 2\nu-1$, then
 \eqref{strong} is true if and only if $\mu\ge {(m-1)/{2}}+\nu+3$.
\item Suppose that for some $n=1,2,3$, we have
  $\varepsilon\ge 2^{1-n}(m+(2\nu-1)(2^{n-1}+1))$,
  ${(m-1)/{2}}+\nu+1-n>0$ and
  $$\mu-n\ge\min\left\{m-1+2\nu+2-2n;\,\max\{1;\frac{m-1}{2}+\nu+1-n\}\right\}.$$
  Then, condition \eqref{strong}  is
 true.
 \end{enumerate}
\end{theorem}
%%%%

We now provide a characterization result for the following problem:
 wether the condition~\eqref{strong} is satisfied for fixed
$\beta_2,\beta_1>0$ (and not for any $\beta_2>\beta_1>0$ as in the
problem~\ref{leo_messi2}).

\begin{theorem}\label{th_fixed}
Let $\mu>0$, $\nu>\frac{1}{2}$, $\mu+\nu>1$ and $m\in\N$. Let
$\varepsilon\in\R$, $\beta_2,\beta_1>0$ and $a:=\beta_2/\beta_1$.
 Then
 $f_{\mu,\nu,\varepsilon,\beta_1,\beta_2}\in\Phi_m$ if and only if
 \begin{equation}\label{fix}
 \frac{1}{x^\mu
(1+x^2)^{\frac{m-1}{2}+\nu} }-
 \frac{a^{2\nu-1-\varepsilon}}{x^\mu
(1+a^2x^2)^{\frac{m-1}{2}+\nu} }\in  {\cal CM}\,.
 \end{equation}
 If $f_{\mu,\nu,\varepsilon,\beta_1,\beta_2}\in\Phi_m$ and $\beta_2>\beta_1>0$, then $\varepsilon\ge2\nu-1$.
\end{theorem}

\noindent {\bf Proof.} From Remarks~\ref{rem_1} and \ref{re5}, it follows
that the following conditions are equivalent:
\begin{enumerate}
\item
$f_{\mu,\nu,\varepsilon,\beta_1,\beta_2}\in\Phi_m$.
\item
  $\beta_2^{\varepsilon+m}\F_m(h_{\mu,\nu})(\beta_2t)-\beta_1^{\varepsilon+m}\F_m(h_{\mu,\nu})(\beta_1t)\ge 0$ for all $t\in(0,\infty)$.
\item
  $L\left(t^{m-1+2\nu+\mu -1}  \left(\beta_2^{\varepsilon+m}\F_m(h_{\mu,\nu})(\beta_2t)-\beta_1^{\varepsilon+m}\F_m(h_{\mu,\nu})(\beta_1t)\right)\right)(x)\in {\cal CM}$.
\end{enumerate}
Furthermore, we have (see~\eqref{pdh}):
\begin{equation*}
 \beta^{\varepsilon+m}L\left(t^{m-1+2\nu+\mu -1} \F_m(h_{\mu,\nu})(\beta t)\right)(x)
 =
 \frac{\beta^{\varepsilon-2\nu+1}}{x^\mu(1+x^2/\beta^2)^{\frac{m-1}{2}+\nu}}
 \;,\; \beta>0,x>0\,.
\end{equation*}
Suppose that condition~\eqref{fix} is satisfied. Since completely monotone functions are non-negative on $(0,+\infty)$,
 we have that $1-a^{2\nu-1-\varepsilon}\ge 0$. If, additionally, $a=\beta_2/\beta_1>1$,
 then $\varepsilon\ge 2\nu-1$. The proof is completed. \hfill
 $\blacksquare$

%%%%%%%%%%%%%%
Direct inspection of the proof of Proposition~\ref{leo_messi3} as well as the
proof of Theorem~\ref{th_fixed}  shows that the following proposition is true.
\begin{prop}\label{eqv}
Let $\mu,\nu,m,2\nu-1+m>0$ and $\varepsilon\in\R$. Then following conditions are equivalent:
 \begin{equation*}\label{eqvfix}
 \begin{split}
 &1.\;\; \frac{\varepsilon-2\nu+1+(\varepsilon+m)x^2}{x^\mu (1+x^2)^{\frac{m-1}{2}+\nu+1}}\in {\cal CM}\,.
 \\
 &2.\;\; \frac{1}{x^\mu
(1+x^2)^{\frac{m-1}{2}+\nu} }-
 \frac{a^{2\nu-1-\varepsilon}}{x^\mu
(1+a^2x^2)^{\frac{m-1}{2}+\nu} }\in  {\cal CM}\,,\forall a>1\,.
\\
 &3.\;\; \frac{1}{x^\mu
(1+x^2)^{\frac{m-1}{2}+\nu} }-
 \frac{a_n^{2\nu-1-\varepsilon}}{x^\mu
(1+a_n^2x^2)^{\frac{m-1}{2}+\nu} }\in  {\cal CM}\,\text{ for some sequence }
a_n>1\,, a_n\to 1\,.
\end{split}
 \end{equation*}
 If the first condition is satisfied, or if the second condition is satisfied for some $a>1$, then $\varepsilon\ge 2\nu-1$.
\end{prop}
{\bf Proof} The crux of the proof is in the following equality, which holds for for $a>1$:
\begin{equation*}\label{eqvfix}
 \int_{1}^{a}\frac{\left(\varepsilon-2\nu+1+(\varepsilon+m\right)x^2t^2)\;t^{2\nu-\varepsilon-2}}{x^\mu
 (1+x^2t^2)^{\frac{m-1}{2}+\nu+1}}\;dt
 =
 \frac{1}{x^\mu
(1+x^2)^{\frac{m-1}{2}+\nu} }-
 \frac{a^{2\nu-1-\varepsilon}}{x^\mu
(1+a^2x^2)^{\frac{m-1}{2}+\nu}}\,.
 \end{equation*}
Since completely monotone functions are closed under scale mixtures, the equality above provides the implication $1\Rightarrow 2$. Assertion $2\Rightarrow 3$ is obvious. To prove assertion
$3\Rightarrow 1$, it is necessary to take in the last equation $a=a_n$, divide both
sides by $a_n-1$ and take the limit as $n\to\infty$, Finally, we make use of the
well-known fact: if a sequence of completely monotone functions converges
pointwise on $(0,+\infty)$, then the limit function is also completely
monotone. \hfill $\blacksquare$

%%%%
We conclude this section detailing the exact smoothness of the differences of Buhmann functions.
\begin{theorem}\label{th_smoothness}
Let $\nu\in\N$, $\mu>0$, $\varepsilon\in\R$, $\beta_2,\beta_1>0$,
and $\beta_2\ne\beta_1$. Let $q:= \min(\beta_1,\beta_2)$. Then,
\begin{enumerate}
\item If $\varepsilon\ne 2\nu-1$, then
$f_{\mu,\nu,\varepsilon,\beta_1,\beta_2}\in C^{2\nu-2}(-q,q)$, and
   $f_{\mu,\nu,\varepsilon,\beta_1,\beta_2}\not\in
   C^{2\nu-1}(-q,q)$.
\item If $\varepsilon= 2\nu-1$, $\mu\notin \{1,2\}$, then
$f_{\mu,\nu,\varepsilon,\beta_1,\beta_2}\in C^{2\nu}(-q,q)$, and
   $f_{\mu,\nu,\varepsilon,\beta_1,\beta_2}\not\in
   C^{2\nu+1}(-q,q)$.
\item  If $\varepsilon= 2\nu-1$, $\mu=1$ or $\mu=2$, then
$f_{\mu,\nu,\varepsilon,\beta_1,\beta_2}$ is
a even polynomial of degree at most $\mu+2\nu-2$ on $[-q,q]$, and therefore
$f_{\mu,\nu,\varepsilon,\beta_1,\beta_2}\in C^{\infty}(-q,q)$.
\end{enumerate}
\end{theorem}

%\proof
\noindent {\bf Proof.}   If $\mu,\nu>0$, then arguments in \cite[Equality~(40)]{c10} show that
\begin{equation}\label{f40}
h_{\mu,\nu}(x)= (1-x)^{\mu+\nu-1}\int_{0}^{1} t^{\mu-1}
(1-t)^{\nu-1} (1-t+(1+t)x)^{\nu-1}\;{\rm d}t\; ,\; x\in (0,1)\,.
\end{equation}
Let $\nu\in\N$. Then from Proposition 1 in \cite{zast2006} we have
that $h_{\mu,\nu}\in C^{2\nu-2}(-1,1)$ and
   $h_{\mu,\nu}\not\in C^{2\nu-1}(-1,1)$. Thus,
   from \eqref{f40} it follows that
   \begin{equation}\label{series}
   \begin{split}
   &
  h_{\mu,\nu}(x)=\sum_{k=0}^{\infty}{a_k(\mu,\nu)}x^{2k}+|x|^{2\nu-1}\sum_{k=0}^{\infty}{b_k(\mu,\nu)
  x^{2k}}\,, \qquad |x|<1,
  \\&
  b_0(\mu,\nu)\ne 0\,,\;\;b_1(\mu,\nu)
  =\frac{h^{(2\nu+1)}_{\mu,\nu}(+0)}{(2\nu+1)!}\,.
  \end{split}
   \end{equation}
From \eqref{rek} we have that $h^{'}_{\mu,\nu}(x)=-2(\nu-1)\,x\,
h_{\mu,\nu-1}(x)$ for $\nu\ge 2$, $x>0$, and
 $h^{(k+1)}_{\mu,\nu}(x)=-2(\nu-1)(x\,h^{(k)}_{\mu,\nu-1}(x)+k\,
h^{(k-1)}_{\mu,\nu-1}(x))$ for $k\ge 1$, $0<x<1$.
 From the last equation it follows that
 $h^{(k+1)}_{\mu,\nu}(+0)=-2(\nu-1)\,k\,
h^{(k-1)}_{\mu,\nu-1}(+0)$ for $\nu\ge 2$, $k\ge 1$, and for $\nu\ge
2$, $k\ge 2\nu-3$, we have (for convenience, we consider that
$0!!:=1$ and $(-1)!!:=1$)
\begin{equation}\label{der}
h^{(k+1)}_{\mu,\nu}(+0)=(-2)^{\nu-1}\,(\nu-1)!\cdot
\frac{k!!}{(k-2\nu+2)!!}\, h^{(k-2\nu+3)}_{\mu,1}(+0)\,.
\end{equation}
  The  equality~\eqref{der} is true  for $\nu=1$, $k\ge-1$.
  It's obvious that $ h_{\mu,1}(x)=(1-|x|)^{\mu}_+/\mu$ and
  $ h^{(p)}_{\mu,1}(+0)=(-1)^p\Gamma(\mu)/\Gamma(\mu-p+1)$, $p\in\Z_+$.
 Therefore, for $k\ge\nu-1$, with $\nu\in\N$, we have
\begin{equation*}
\begin{split}
h^{(2k+1)}_{\mu,\nu}(+0)
&=(-4)^{\nu-1}(\nu-1)!\cdot\frac{k!}{(k-\nu+1)!}\cdot
h^{(2k-2\nu+3)}_{\mu,1}(+0)
\\&=
 -(-4)^{\nu-1}(\nu-1)!\cdot\frac{k!}{(k-\nu+1)!}\cdot\frac{\Gamma(\mu)}{\Gamma(\mu-2k+2\nu-2)}\,;
\\
 h^{(2\nu+1)}_{\mu,\nu}(+0)&=-(-4)^{\nu-1}(\nu-1)!\nu!(\mu-1)(\mu-2)\,.
\end{split}
\end{equation*}
On the other hand, \eqref{series} shows that, for
$|x|<q=\min(\beta_1,\beta_2)$,
\begin{eqnarray*}
 f_{\mu,\nu,\varepsilon,\beta_1,\beta_2}(x)
& =& \sum_{k=0}^{\infty}{a_k(\mu,\nu)}\left(\beta_2^{\varepsilon-2k}-\beta_1^{\varepsilon-2k}\right)x^{2k} + \\ &+& |x|^{2\nu-1}\sum_{k=0}^{\infty}{b_k(\mu,\nu)
  \left(\beta_2^{\varepsilon-2\nu+1-2k}-\beta_1^{\varepsilon-2\nu+1-2k}\right)x^{2k}}\,.
\end{eqnarray*}
Assertions 1. and 2. are thus proved.

If, in addition, $\mu\in\N$, then $h_{\mu,\nu}$
 is a polynomial of degree $\mu+2\nu-2$ on the compact interval $[0,1]$. If $\mu=1$ or
$\mu=2$, then in the second sum of \eqref{series}, the terms with
$k\ge 1$ vanish. Therefore, if $\varepsilon= 2\nu-1$, $\mu=1$ or
$\mu=2$, then $f_{\mu,\nu,\varepsilon,\beta_1,\beta_2}$ on the
interval $[-q,q]$ is a even polynomial of degree $\le\mu+2\nu-2$.
Assertion 3.  is proved. The proof is completed.  $\blacksquare$
%\endproof

\section{Connections with previous literature}

 The $k$-fold application of the Mont{\'e}e operator to the Askey function ${\cal A}_{\mu}(t)= (1-t)_+^{\mu}$, for $\mu \ge (m+1)/2$, results in Wendland functions $\psi_{\mu,k}$ as described in previous section.
Table 1 depicts the role of the mapping $f_{\mu,\nu,\varepsilon,\beta_1,\beta_2}$ as defined through Equation (\ref{del_piero}). In particular, we consider $f_{\mu,k+1,2k+1,\beta_1,\beta_2}$, for $k=0,1,2$. Expressions of the corresponding Wendland functions are reported in the second column of the same table. According to Theorem \ref{cr7}, $f_{\mu,k+1,2k+1,\beta_1,\beta_2}\in \Phi_m$ if and only if $\mu \ge (m+7)/2+k$. The third and fourth column allow to describe the action of the mapping $f_{\mu,k+1,2k+1,\beta_1,\beta_2}$. When $\mu \notin \{1,2\}$, one can clearly appreciate the increase in terms of differentiability at the origin.

\begin{table}[!b] \label{francesco_totti}
\begin{center}
\begin{tabular}{|c|c|c|c|}
\hline
$k$ & $\psi_{\mu,k}$ & {D}$_{{\tt before}}$ & {D}$_{{\tt after}}$ \\
& & & \\
\hline
&&& \\
0 & $(1-x)_{+}^{\mu}$ & $C(\{0\})$ & $C^{{2}}(\{0\})$ \\
&&& \\
1 & $(1-x)_+^{\mu+1}(1+(\mu+1)x)$ & $C^{2}(\{0\})$ & $C^{{4}}(\{0\})$ \\
&&& \\
$2 $ & $(1-x)_+^{\mu+2}(1+(\mu+2)x+ \frac{1}{3}\left ( (\mu+2)^2-1\right ))$ & $C^{4}(\{0\})$ & $C^{{6}}(\{0\})$  \\
&&& \\
\hline
\end{tabular}
\vspace{0.2cm}
\caption{Examples of the Wendland functions $\psi_{\mu,k}$ for $k=0,1,2$ and $\mu \notin \{1,2 \}$. For all cases, $x \ge 0$. The first column reports the values of $k$, and the corresponding expression in the second column reports the analytic expression of $\psi_{\mu,k}$. In the third column, {\em D}$_{{\tt before}}$ depicts the differentiability of $\psi_{\mu,k}$. In the last column, {\em D}$_{{\tt after}}$ stays for the differentiability at the origin of $f_{\mu,k+1,2k+1,\beta_1,\beta_2}$. Observe that when $\mu=1$ or $\mu=2$, then {\em D}$_{{\tt after}}=\infty$.}
\end{center}
\end{table}

\begin{figure}
\begin{tabular}{ccc}
  \includegraphics[scale=0.30]{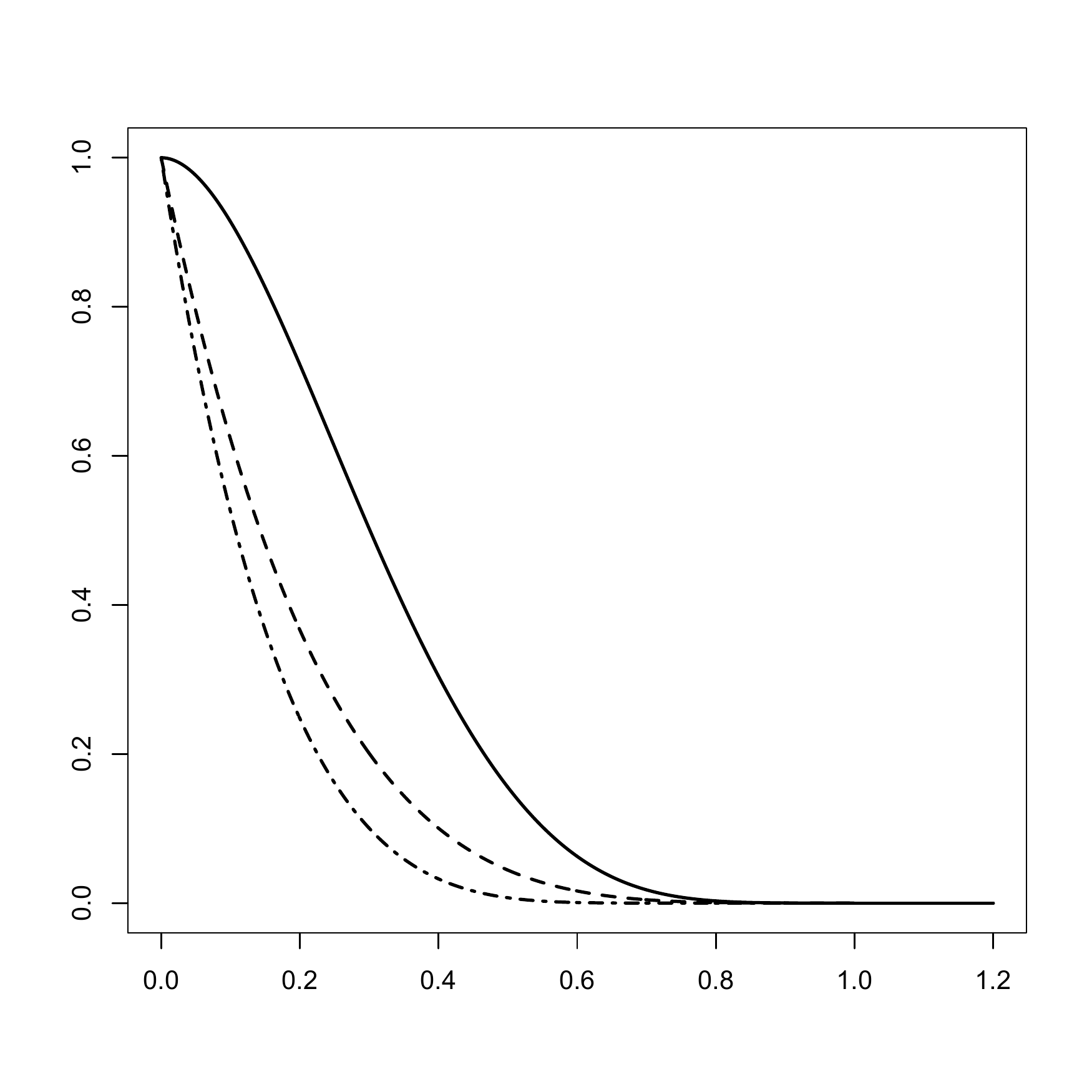} & \includegraphics[scale=0.30]{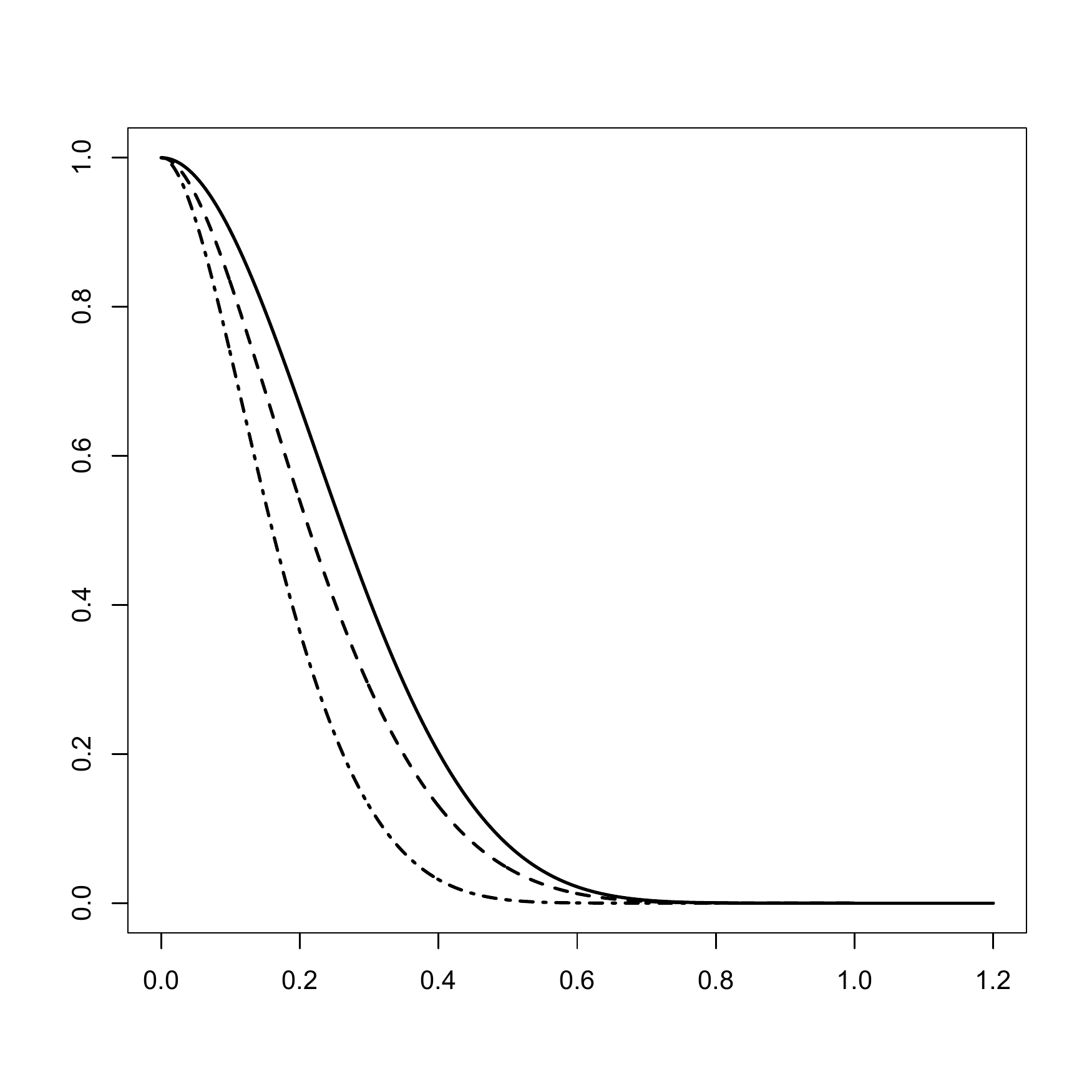} & \includegraphics[scale=0.30]{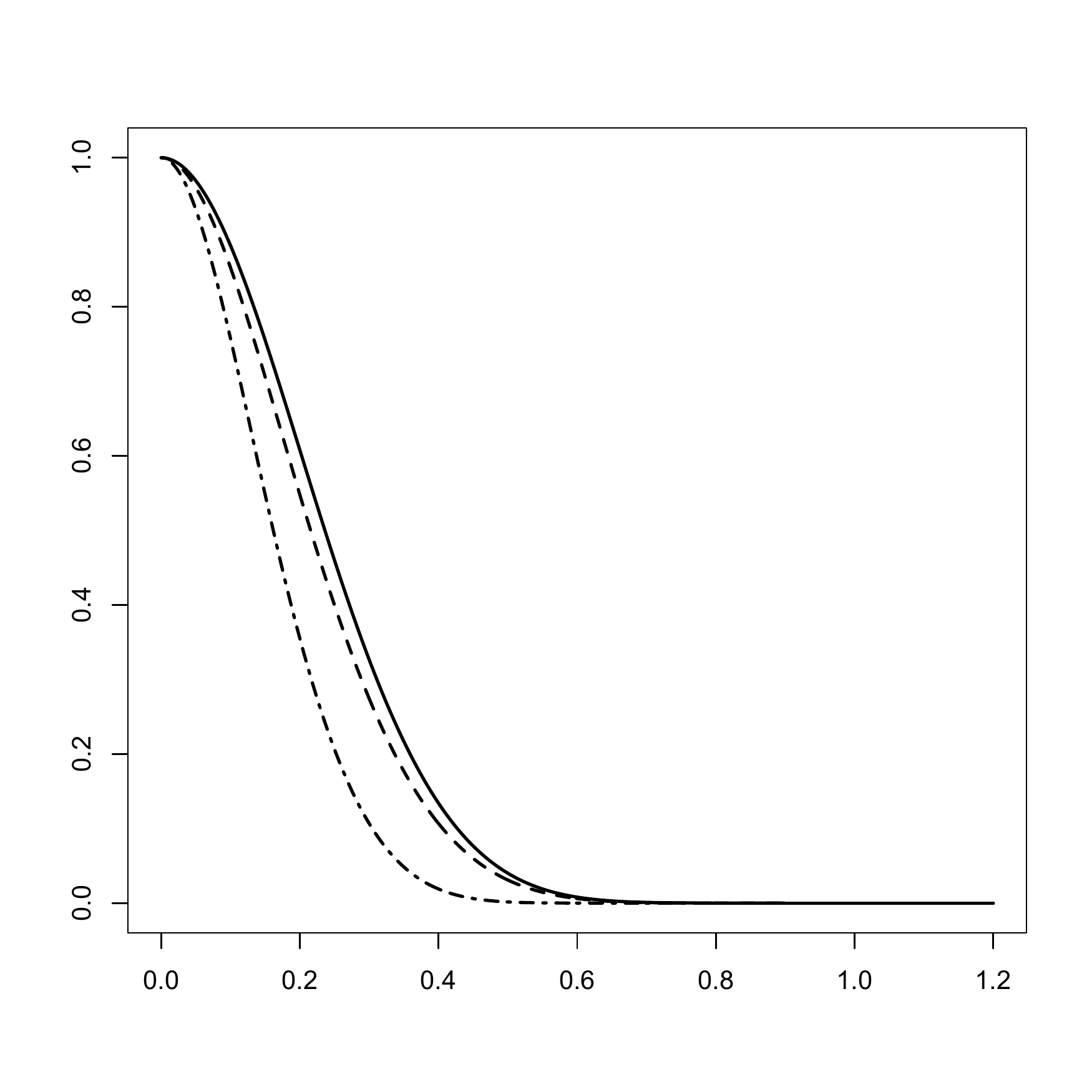}
\end{tabular}
\caption{Continuous Lines: $f_{(d+1)/2+k+3,k+1,2k+1,0.75,1}$, defined according to Equation (\ref{del_piero}) for $k=0, 1, 2$ (from left to right). Dashed and Dashed-dot Lines report Wendland functions with $k=0,1,2$ and compact supports $\beta_1=0.75$ and $\beta_2=1$. All the functions are normalized with their value at the origin.}
\end{figure}

The Wendland radial basis functions are piecewise polynomial compactly supported reproducing kernels in Hilbert spaces which are norm–equivalent to Sobolev spaces. But they only cover the Sobolev spaces
$H^{d/2+k+1/2}(\R^d)$, when $k \in \mathbb{N}$. Motivated by this fact, Robert Schaback \cite{Schaback} covered the case of the integer order spaces in even dimensions. Namely, he derived the missing Wendland functions working for half-integer $k$ and even dimensions, reproducing integer-order Sobolev spaces in even dimensions, and showing that they turn out to have two additional non-polynomial terms: a logarithm and a square root.

Other walks through dimensions have been recently proposed by \cite{porcu-zast} through the  Generalized Askey functions
$\varphi_{n,k,m}:[0,\infty) \R$ defined through
$$ \varphi_{n,k,m}(t)= t^{k-n}(1-t)_+^{n+m+1} F(
n - k, n + 1, n + m + 2, 1 -1/t), \qquad t \ge 0, $$ with $F$ being the Gauss hypergeometric function. The parameters $(n,m,k)$ are then shown to be crucial in order to determine when $\Phi_{m}$ (see their Proposition 2.3).

To our knowledge, the only case of compactly supported correlation functions which is not covered by this work is the case of the Euclid's hat (\cite{euclid1, euclid2}), which is the self-convolution of an indicator function supported on the unit ball in $\R^d$. As noted by \cite{Schaback}, while Euclid’s hat is not differentiable and Wu’s functions have zeros in their Fourier transform, Wendland’s functions have no such drawbacks. They are polynomials on $[0, 1]$ and yield positive definite $2k$-times differentiable  radial basis functions on $\R^d$. Given these properties, their polynomial degree $\lfloor d/2\rfloor + 3k + 1$ is minimal.

The proof of Theorem \ref{cr7} highlights explicit connections with previous literature devoted to (sub) classes of completely monotone functions. A function $f:[0,\infty) \to \R$ is called Logarithmically completely monotonic on $(0,\infty)$, and denoted $f \in {\cal L}(0,\infty)$ if and only if it is infinitely often differentiable on $(0,\infty)$ and
$$ (-1)^n \left [ \log f(x)\right]^{(n)} \ge 0, \qquad x \ge 0. $$
By well known results (see \cite{bergPM}, with the references therein) $f \in {\cal L} \Leftrightarrow f^{\alpha} \in {\cal CM} \Leftrightarrow (f)^{1/n} \in {\cal CM}$, for all $\alpha >0$ and $n \in \mathbb{N}$.

Berg, Porcu and Mateu \cite{bergPM} introduced the so called {\em Auxiliary} family
$$ f_{\alpha,\beta}(x) = \frac{1}{x^{\alpha}(1+x^{\beta})}, \qquad x >0, $$
with  $\alpha,\beta$ positive parameters. They show that
$f_{\alpha,\beta} \in {\cal L}$ for $\alpha \ge 0$ and $0 \le \beta
\le 1$. Additionally, $f_{\alpha,2} \in {\cal L}$ if and only if
$\alpha \ge 2$. This resut is one of the crux for describing the
complete monotonicity of the Dagum family. Known cases, when
$x^{-\mu}(1+x^2)^{-\nu}\notin {\cal CM}$ (if $\mu,\nu>0$, then this
is equivalent to the inequality $I_{1,\mu,\nu,2\nu-1}\ge 0$ is not
true for all $t>0$): 1) $\mu<0$ or $\mu=0$, $\nu\ne 0$ (it is
obvious); 2) $\mu<\nu$; 3) $0<\mu=\nu<1$. Proof of the latter two
cases, see, for example,  \cite{moak}, \cite[Lemma~8]{Zast2000} and
\cite[Theorem 5]{zast2006}.

\section*{Acknowledgement}
We are grateful to Professor R. Furrer for helpful remarks through the preparation of the manuscript.
Research work of Moreno Bevilacqua was partially supported by grant FONDECYT  11121408 from Chilean government.
Research work of Emilio Porcu was partially supported by grant  FONDECYT 1130647 from Chilean government.

{%\small

}
%%%%%%%%

\end{document}